\documentclass[12pt, a4paper]{amsart}
\usepackage{amsmath}
\usepackage{geometry,amsthm,graphics,tabularx,amssymb,shapepar}
\usepackage{amscd}
\usepackage[all,2cell,dvips]{xy}

\newcommand{\CN}{{\mathcal {N}}}

\newcommand{\CU}{{\mathcal {U}}}

\newcommand{\End}{{\mathrm{End}}}

\newcommand{\GL}{{\mathrm{GL}}}

\newcommand{\Hom}{{\mathrm{Hom}}}

\newcommand{\rank}{{\mathrm{rank}}}

\newcommand{\SO}{{\mathrm{SO}}}

\newcommand{\vsp}{{\vspace{0.2in}}}

\newcommand{\con}{\textit{C}}

\newcommand{\rs}{\operatorname{s}}

\newcommand{\oO}{\operatorname{O}}
\newcommand{\oS}{\operatorname{S}}

\newcommand{\oU}{\operatorname{U}}
\newcommand{\oZ}{\operatorname{Z}}

\renewcommand{\u}{\mathfrak u}

\newcommand{\z}{\mathfrak z}

\newcommand{\C}{\mathbb{C}}
\newcommand{\R}{\mathbb R}

\newcommand{\la}{\langle}
\newcommand{\ra}{\rangle}

\newcommand{\be}{\begin {equation}}
\newcommand{\ee}{\end {equation}}
\newcommand{\bee}{\begin {equation*}}
\newcommand{\eee}{\end {equation*}}

\theoremstyle{Theorem}

\newtheorem{thm}{Theorem}[section]

\newtheorem{lemt}[thm]{Lemma}
\newtheorem{prpt}[thm]{Proposition}

\theoremstyle{Theorem}

\theoremstyle{Theorem}

\theoremstyle{Plain}

\theoremstyle{Definition}

\begin{document}

\title[Special orthogonal groups]{A note on special orthogonal groups following Waldspurger}

\author{Binyong Sun}
\address{Academy of Mathematics and Systems Science\\
Chinese Academy of Sciences\\
Beijing, 100190, P.R. China} \email{sun@math.ac.cn}

\author{Chen-Bo Zhu}
\address{Department of Mathematics\\
National University of Singapore\\
2 Science drive 2\\
Singapore 117543} \email{matzhucb@nus.edu.sg}

\thanks{First version on May 14, 2010}

\begin{abstract}
The purpose of this note is to verify that the archimedean multiplicity one theorems shown for orthogonal groups (as well as general linear and unitary groups) in a previous paper of the authors remain valid for special orthogonal groups. The necessary ingredients to establish this variant are due to Waldspurger.
\end{abstract}

\maketitle

\begin{thm}\label{main}
Let $G$ be a special orthogonal group $\SO(p,q)$ or $\SO_n(\C)$, $p,n\geq 1$, $q\geq 0$. Let $G'$ be $\SO(p-1,q)$ or $\SO_{n-1}(\C)$, viewed as a subgroup of $G$ as usual. Then for every irreducible Casselman-Wallach smooth representation $V$ of $G$, and $V'$ of $G'$, one has that
\[
  \dim \Hom_{G'}(V\widehat \otimes V',\C)\leq 1.
  \]
  Here ``$\widehat \otimes$" stands for the completed projective tensor product of Hausdorff locally convex topological vector spaces.
\end{thm}

\vsp

We follow the general set-up of \cite[Section 3]{SZ}.

Let $(A,\tau)$ be a (finite-dimensional) commutative involutive algebra over $\R$, and let $E$ be a (non-degenerate finitely generated) Hermitian $A$-module, with a Hermitian form
\[
  \la \,,\,\ra_E:E\times E\rightarrow A.
\]
Denote by $\oU(E)$ the group of $A$-linear automorphisms of $E$ preserving the form $\la\,,\,\ra_E$.  Write $E_\R:=E$, viewed as a real vector space. Denote by $\breve \oU(E)$ the subgroup of $\GL(E_\R)\times \{\pm 1\}$ consisting of pairs $(g,\delta)$ such that either
\[
  \delta=1 \quad \textrm{and}\quad \la gu,gv\ra_E=\la u,v\ra_E, \quad u,v\in E,
\]
or
\[
  \delta=-1 \quad \textrm{and}\quad \la gu,gv\ra_E=\la v,u\ra_E, \quad u,v\in E.
\]
This contains $\oU(E)$ as a subgroup of index two.

\vsp
First assume that $(A,\tau)$ is simple. If $\tau$ is nontrivial, we put
\[
  \oU_{\rs}(E):=\oU(E) \quad \textrm{(This is a general linear group or a unitary group.)}\quad
\]
and
\[
 \breve  \oU'_{\rs}(E):=\breve \oU_{\rs}(E):=\breve \oU(E).
 \]
Otherwise, $\tau$ is trivial and $A=\R$ or $\C$. Then we have $\oU(E)=\oO(E)$, and $\breve \oU(E)=\oO(E)\times \{\pm 1\}$, in the usual notations.
We shall put
\[
  \oU_{\rs}(E):=\oS \oO(E)\subset \oO(E),
\]
and following Waldspurger \cite{Wa}
\[
   \breve \oU_{\rs}(E):=\left \{(g,\delta)\in \breve \oU(E)=\oO(E)\times \{\pm 1\}\mid \det(g)=\delta^{\left[\frac{\dim_A E +1}{2} \right]}\right \}
\]
and
\[
   \breve \oU'_{\rs}(E):=\left \{(g,\delta)\in \breve \oU(E)=\oO(E)\times \{\pm 1\}\mid \det(g)=\delta^{\left[\frac{\dim_A E}{2} \right]}\right \}.
\]
In general, write
\[
  (A,\tau)=(A_1,\tau_1)\times (A_2,\tau_2)\times \cdots\times (A_r,\tau_r)
\]
as a product of simple commutative involutive algebras over $\R$. Then
\begin{equation}\label{Eproduct}
  E=E_1\times E_2\times \cdots\times E_r,
\end{equation}
where
\[
  E_i:=A_i\otimes_A E
\]
is naturally a Hermitian $A_i$-module. We put
\[
  \oU_{\rs}(E):=\oU_{\rs}(E_1)\times \oU_{\rs}(E_2)\times \cdots \times \oU_{\rs}(E_r)\subset \oU(E),
\]
and
\begin{eqnarray*}
       \breve \oU_{\rs}(E)&:=&\breve \oU_{\rs}(E_1)\times_{\{\pm 1\}} \breve \oU_{\rs}(E_2)\times_{\{\pm 1\}} \cdots \times _{\{\pm 1\}} \breve \oU_{\rs}(E_r)\\
       &:=&\{(g_1,g_2,\cdots,g_r,\delta)\mid (g_i,\delta)\in \breve \oU_{\rs}(E_i),\, i=1,2,\cdots, r \}\\
       &\subset & \breve \oU(E).
\end{eqnarray*}
The latter ($\breve \oU_{\rs}(E)$) contains the former as a subgroup of index two. Denote by $\chi_{\rs, E}$ the quadratic character on $\breve \oU_{\rs}(E)$ with kernel $\oU_{\rs}(E)$. Likewise we define a group $\breve \oU'_{\rs}(E)$ which contains $\oU_{\rs}(E)$ as a subgroup of index two. Denote by $\chi'_{\rs, E}$ the quadratic character on $\breve \oU'_{\rs}(E)$ with kernel $\oU_{\rs}(E)$.

\vsp

Write
\[
  \u_{\rs}(E):=\{x\in \End_A(E)\mid \la xu,v\ra_E+\la u, xv\ra_E=0,\,u,v\in E\}
\]
for the Lie algebra of $\oU_{\rs}(E)$ (which is also the Lie algebra of $\oU(E)$). Let the groups $\breve \oU_{\rs}(E)$ and $\breve \oU'_{\rs}(E)$ act on $\oU_{\rs}(E)$ and $\u_{\rs}(E)$ by
\begin{equation}\label{actioninf}
   \left\{
     \begin{array}{ll}
        (g,\delta)\cdot x:=gx^{\delta}g^{-1},\quad& x\in \oU_{\rs}(E),\medskip\\
        (g,\delta)\cdot x:=\delta \,gxg^{-1},\quad& x\in \u_{\rs}(E).
    \end{array}
  \right.
 \end{equation}
Also they act on $E$ by
\begin{equation}\label{actioninf}
   \left\{
     \begin{array}{ll}
        (g,\delta)\cdot u:=\delta \,gu,\quad & (g,\delta)\in \breve \oU_{\rs}(E),\,u\in E,\medskip\\
        (g,\delta)\cdot u:=gu,\quad & (g,\delta)\in \breve \oU'_{\rs}(E),\, u\in E.
    \end{array}
  \right.
 \end{equation}

It is by now standard (see for example \cite[Section 7]{SZ}) that Theorem \ref{main} is implied by the first assertion of the following theorem in the case of $A=\R$ or $\C$, $\tau$ trivial.

\begin{thm}\label{vanish1}
One has that
\begin{equation}\label{ev1}
  \con^{-\infty}_{\chi_{\rs,E}}(\oU_{\rs}(E)\times E)=0
\end{equation}
and
\begin{equation}\label{ev2}
  \con^{-\infty}_{\chi_{\rs,E}}(\u_{\rs}(E)\times E)=0.
\end{equation}
\end{thm}

Note that $(g,\delta)\mapsto (\delta g,\delta)$ is a group isomorphism from $\breve \oU_{\rs}(E)$ onto $\breve \oU'_{\rs}(E)$ fixing $\oU_{\rs}(E)$. Thus Theorem \ref{vanish1} is equivalent to
 \begin{thm}\label{vanish2}
One has that
\[
  \con^{-\infty}_{\chi'_{\rs,E}}(\oU_{\rs}(E)\times E)=0
\]
and
\[
  \con^{-\infty}_{\chi'_{\rs,E}}(\u_{\rs}(E)\times E)=0.
\]
\end{thm}

\vsp

For $E$ as in \eqref{Eproduct}, put
\[
  \operatorname{sdim}(E):=\sum_{i=1}^r \max \{\rank_{A_i} E_i -1, 0\}+\dim_\R E_\R.
\]
We argue by induction on $\operatorname{sdim}(E)$ and so will assume that Theorem \ref{vanish1} (and hence Theorem \ref{vanish2}) holds whenever
$\operatorname{sdim}(E)$ is smaller.

Without loss of generality, in the remaining part of this note, assume that $(A,\tau)$ is simple and $E$ is faithful as an $A$-module. Let $x$ be a semisimple element of $\oU_{\rs}(E)$ or $\u_{\rs}(E)$. Denote by $A_x$ the subalgebra of $\End_A(E)$ generated by $A$, $x$ and $x^\tau$. Here $\tau$ is the involution of $\End_A(E)$ given by
\[
  \la xu,v\ra_E=\la u, x^\tau v\ra_E,\quad u,v\in E.
\]
Then $A_x$ is again a commutative involutive algebra over $\R$, and $E_x:=E$ is naturally a Hermitian $A_x$-module.
In the notations of this note, the following lemma is the first key observation of Waldspurger \cite{Wa}.

\begin{lemt}
The group $\oU_{\rs}(E_x)$ is a subgroup of $\oU_{\rs}(E)$, the Lie algebra $\u_{\rs}(E_x)$ is a Lie subalgebra of $\u_{\rs}(E)$, and $\breve \oU_{\rs}(E_x)$ is a subgroup of $\breve \oU_{\rs}(E)$.  The embeddings $\oU_{\rs}(E_x)\hookrightarrow \oU_{\rs}(E)$ and $\u_{\rs}(E_x)\hookrightarrow \u_{\rs}(E)$ are both  $\breve \oU_{\rs}(E_x)$-equivariant. Furthermore $x\in \oU_{\rs}(E_x)$ if $x\in \oU(E)$ and $x\in \u_{\rs}(E_x)$ if $x\in \u(E)$.
\end{lemt}

As in \cite[Section 5]{SZ}, Harish-Chandra's method of descent and the above lemma imply the following

\begin{prpt}\label{reduc1}
Every element of $\con^{-\infty}_{\chi_{\rs,E}}(\oU_{\rs}(E)\times E)$ is supported in $(\oZ_E\times \CU_E)\times E$,
and every element of $\con^{-\infty}_{\chi_{\rs,E}}(\u_{\rs}(E)\times E)$ is supported in $(\z_E\oplus \CN_E)\times E$, where $\oZ_E$ is the scalar multiplications (by $A$) in $\oU_{\rs}(E)$, $\z_E$ is the scalar multiplications (by $A$) in $\u_{\rs}(E)$, $\CU_E$ is the set of unipotent elements of $\oU_{\rs}(E)$, and $\CN_E$ is the set of nilpotent elements of $\u_{\rs}(E)$.
\end{prpt}

By the first assertion of the above proposition, \eqref{ev2} will imply \eqref{ev1}. So we only need to prove \eqref{ev2}.

\vsp

Let $v$ be a non-degenerate element of $E$ (i.e., $\la v,v\ra_E$ is invertible in $A$), and denote by $E_v$ the orthogonal complement of $v$ in $E$. The second key observation of Waldspurger \cite{Wa} is the following
\begin{lemt}
The map $(g,\delta)\mapsto (g|_{(E_v)_\R}, \delta)$ identifies the stabilizer of $v$ in $\breve \oU_{\rs}(E)$ with the group $\breve \oU'_{\rs}(E_v)$. Furthermore, the restriction to $\breve \oU'_{\rs}(E_v)$ of the module $\u_{\rs}(E)$ is isomorphic to $\u_{\rs}(E_v)\times E_v\times \u_{\rs}(Av)$. Here $\u_{\rs}(Av)$ carries the trivial $\breve \oU'_{\rs}(E_v)$-action.
\end{lemt}

Again Harish-Chandra's method of descent and the above lemma imply the following

\begin{prpt}\label{reduc2}
Every element of $\con^{-\infty}_{\chi_{\rs,E}}(\u_{\rs}(E)\times E)$ is supported in $\u_{\rs}(E)\times \Gamma_E$, where $\Gamma_E:=\{u\in E\mid \la u,u\ra_E=0\}$ is the null cone of $E$.
\end{prpt}

The (same and key) argument of \cite[Section 4]{SZ} (reduction within the null cone) works in the setting of this note and we have
\begin{prpt}\label{reducn}
Assume that every element of $\con^{-\xi}_{\chi_{\rs,E}}(\u_{\rs}(E)\times E)$ is supported in $(\z_E\oplus \CN_E)\times \Gamma_E$, then
\begin{equation}\label{evt}
  \con^{-\xi}_{\chi_{\rs,E}}(\u_{\rs}(E)\times E)=0.
\end{equation}
Here $\con^{-\xi}_{\chi_{\rs,E}}(\u_{\rs}(E)\times E)$ denotes the subspace of tempered generalized functions in $\con^{-\infty}_{\chi_{\rs,E}}(\u_{\rs}(E)\times E)$.
\end{prpt}

Now Propositions \ref{reduc1} and \ref{reduc2} imply that the hypothesis of Proposition \ref{reducn} is satisfied. (This completes the step of reduction to the null cone.) Together with Proposition \ref{reducn}, they imply that \eqref{evt} always holds. Then a general principle due to Aizenbud and Gourevitch (\cite[Theorem 4.0.2]{AGS}) implies that \eqref{ev2} also holds. \qed

\end{document}